December 08, 2016

# Lah numbers and Laguerre polynomials of order negative one

**Khristo N. Boyadzhiev**
Department of Mathematics and Statistics,
Ohio Northern University, Ada, OH 45810, USA
k-boyadzhiev@onu.edu

**Abstract**  In this article we point out interesting connections among Lah numbers, Laguerre polynomials of order negative one, and exponential polynomials. We also discuss several different expressions for the higher order derivatives of *exp* (1/*x*). A new representation of these derivatives is given in terms of exponential polynomials.

## 1. Introduction

The Lah numbers $L(n,k)$ (named after Ivo Lah, a Slovenian mathematician) can be defined by the formula

$$(1.1) \quad L(n,k) = \frac{n!}{k!}\binom{n-1}{k-1}, \quad 1 \leq k \leq n,\ L(0,0) = 1,$$

or, by the generating function

$$\frac{1}{k!}\left(\frac{t}{1-t}\right)^k = \sum_{n=k}^{\infty} L(n,k) \frac{t^n}{n!}.$$

The Lah numbers convert the falling factorial to the rising factorial and vise-versa

$$x(x+1)\ldots(x+n-1) = \sum_{k=1}^{n} L(n,k)\, x(x-1)\ldots(x-k+1),$$

$$x(x-1)\ldots(x-n+1) = \sum_{k=1}^{n}(-1)^{n-k} L(n,k)\, x(x+1)\ldots(x+k-1)$$

(these are the fundamental identities obtained by Ivo Lah).

The Lah numbers have many other interesting applications in analysis and combinatorics (see [1], [2], [9], [12], [16]). They have appeared recently in several papers concerning the consecutive derivatives of the function $\exp(1/x)$. In [10] five proofs were given of the following formula:

(1.2) $\quad D^n e^{1/x} = (-1)^n e^{1/x} x^{-n} \sum_{k=1}^{n} L(n,k)\, x^{-k}$ .

where $D = \dfrac{d}{dx}$ and $n \geq 1$. This formula was obtained also by Feng Qi (see [13] and the remarks in Section 5 there). The formula is a nice application of Lah numbers to a problem in analysis.

At the same time, entry 1.1.3.2 on p.4 in Brychkov's handbook [6] says that

$$\frac{d^n}{dx^n}[x^\lambda e^{-a/x}] = (-1)^n n!\, e^{-a/x} x^{\lambda-n} L_n^{(-\lambda-1)}(a/x)$$

where $L_n^{(\alpha)}(x)$ are the generalized (or associated) Laguerre polynomials of order $\alpha$ (see [14], [16]). The same formula appears as entry 18.5.6. on page 446 in the handbook [15]. With $\lambda = 0$ and $a = -1$ this becomes

(1.3) $\quad D^n e^{1/x} = (-1)^n n!\, e^{1/x} x^{-n} L_n^{(-1)}(-1/x)$ .

As a matter of fact, the derivatives $D^n e^{1/x}$ have been evaluated long time ago. For example, even more general derivatives can be found in the nice little book of Schwatt [18], first published in 1924. The formula on top of page 22 in [18] reads

(1.4) $\quad D^n e^{cx^p} = n!\, e^{cx^p} x^{-n} \sum_{k=1}^{n} \dfrac{(-1)^k}{k!} c^k x^{pk} \sum_{j=1}^{k}(-1)^j \binom{k}{j}\binom{pj}{n}$

where $c, p$ are arbitrary parameters. The same formula appears also on page 27 of that book. With $c = 1$ and $p = -1$ this becomes

$$(1.5) \quad D^n e^{1/x} = (-1)^n n! e^{1/x} x^{-n} \sum_{k=1}^{n} \frac{(-1)^k}{k!} x^{-k} \left\{ \sum_{j=1}^{k} (-1)^j \binom{k}{j} \binom{n+j-1}{n} \right\}$$

since $\binom{-j}{n} = (-1)^n \binom{n+j-1}{n}$.

In the next three sections we discuss the relations among the three formulas for $D^n e^{1/x}$, namely, equations (1.2), (1.3), and (1.5). This will reveal interesting connections of Lah numbers to Laguerre polynomials and also to Stirling numbers. In Section 4 we present a new formula for $D^n e^{cx^p}$ in terms of the exponential polynomials $\varphi_n(x)$ considered in [4] and [5].

## 2. Laguerre polynomials

The generalized Laguerre polynomials can be defined by the generating function

$$\frac{1}{(1-t)^{\alpha+1}} e^{\frac{-xt}{1-t}} = \sum_{n=0}^{\infty} L_n^{(\alpha)}(x) t^n, |t| < 1 ,$$

or by the Rodriguez formula

$$L_n^{(\alpha)}(x) = \frac{e^x x^{-\alpha}}{n!} D^n (e^{-x} x^{n+\alpha}) = \frac{x^{-\alpha}}{n!} (D-1)^n x^{n+\alpha}, \quad n = 0, 1, \ldots$$

(see [14]). When $\alpha = 0$ these are the usual Laguerre polynomials $L_n^{(0)}(x) = L_n(x)$. Usually, in the theory of $L_n^{(\alpha)}(x)$ the restriction $\operatorname{Re} \alpha > -1$ is imposed. In fact, the case $\alpha = -1$ is very interesting and most of the theory holds true for $\alpha = -1$. We shall focus here on the polynomials $L_n^{(-1)}(x)$ defined by

$$L_n^{(-1)}(x) = \frac{xe^x}{n!} D^n (e^{-x} x^{n-1}) = \frac{x}{n!} (D-1)^n x^{n-1}, n = 0, 1, \ldots$$

or, by the generating function, $|t|<1$

$$(2.1) \quad e^{\frac{-xt}{1-t}} = \sum_{n=0}^{\infty} L_n^{(-1)}(x) t^n .$$

We have

$$L_0^{(-1)}(x) = 1,$$

$$L_1^{(-1)}(x) = -x,$$

$$L_2^{(-1)}(x) = \frac{x^2}{2} - x,$$

$$L_3^{(-1)}(x) = \frac{-x^3}{6} + x^2 - x,$$

etc. The coefficients of these polynomials are very close to the Lah number and we can see the exact connection when we compare (1.2) to (1.3). We shall give an independent proof of this connection in order to justify the value $\alpha = -1$ in (1.3).

**Proposition 1**. For any $n \geq 0$,

$$(2.2) \quad L_n^{(-1)}(x) = \frac{1}{n!} \sum_{k=0}^{n} L(n,k) (-x)^k .$$

This reveals the connection between the Lah numbers and the Laguerre polynomials $L_n^{(-1)}(x)$ and it becomes clear now that formulas (1.2) and (1.3) are the same. We also notice that (1.2) is true also for $n=0$ with the summation starting from $k=0$, that is,

$$D^n e^{1/x} = (-1)^n e^{1/x} x^{-n} \sum_{k=0}^{n} L(n,k) x^{-k} .$$

*Proof.* From the Rodriguez formula for $L_n^{(\alpha)}(x)$ one derives easily the representation

$$L_n^{(\alpha)}(x) = \Gamma(n+\alpha+1) \sum_{k=0}^{n} \frac{(-x)^k}{\Gamma(k+\alpha+1) k! (n-k)!}$$

where we cannot set $\alpha = -1$ directly. However, when $n = 0$ this becomes

$$L_0^{(\alpha)}(x) = \frac{\Gamma(\alpha+1)}{\Gamma(\alpha+1)} = 1$$

and any restriction on $\alpha$ can be dropped. For $n \geq 1$ we separate the first term with $k = 0$ and write

$$L_n^{(\alpha)}(x) = \frac{\Gamma(n+\alpha+1)}{\Gamma(\alpha+1)} + \Gamma(n+\alpha+1) \sum_{k=1}^{n} \frac{(-x)^k}{\Gamma(k+\alpha+1) k!(n-k)!} .$$

Setting $\alpha \to -1$ we find for $n \geq 1$

$$L_n^{(-1)}(x) = \Gamma(n) \sum_{k=1}^{n} \frac{(-x)^k}{\Gamma(k) k!(n-k)!} ,$$

since

$$\lim_{\alpha \to -1} \frac{1}{\Gamma(\alpha+1)} = 0 .$$

This representation can be written in the form

(2.3) $\quad L_n^{(-1)}(x) = \sum_{k=1}^{n} \binom{n-1}{k-1} \frac{(-x)^k}{k!}$

which is (2.2). The proof is completed.

The representation (2.2) also shows an important difference between $L_n^{(-1)}(x)$ and $L_n^{(\alpha)}(x)$ for $n \geq 1$. While

$$L_n^{(\alpha)}(0) = \frac{\Gamma(n+\alpha+1)}{\Gamma(\alpha+1)}$$

is different from zero when $\alpha \neq -1$, we have $L_n^{(-1)}(0) = 0$. At the same time, many properties of $L_n^{(\alpha)}(x)$ are shared also by $L_n^{(-1)}(x)$. For example, we have the orthogonally relation ([14, p.84], [17, p. 204-205])

$$\int_0^\infty x^\alpha e^{-x} L_n^{(\alpha)}(x) L_m^{(\alpha)}(x)\, dx = \frac{\Gamma(n+\alpha+1)}{n!}\, \delta_{n,m}$$

for all $n,m \geq 0$ and $\alpha > -1$. Analyzing the proof of this equation in [17] we conclude that it extends to $\alpha = -1$ when $n,m \geq 1$,

(2.4) $$\int_0^\infty x^{-1} e^{-x} L_n^{(-1)}(x) L_m^{(-1)}(x)\, dx = \frac{\delta_{n,m}}{n}\ .$$

This and other properties of $L_n^{(-1)}(x)$ can be used to derive properties for the Lah numbers. Here we have the following:

**Proposition 2**. For any integers $n,m \geq 1,\ n \neq m$

(2.5) $$\sum_{k=1}^n (-1)^k L(n,k) \sum_{j=1}^m (-1)^j L(m,j)(k+j-1)! = 0$$

and when $m = n$

(2.6) $$\sum_{k=1}^n \sum_{j=1}^n (-1)^{k+j} L(n,k)\, L(n,j)\, (k+j-1)! = \frac{(n!)^2}{n}\ .$$

*Proof.* Substituting (2.2) in (2.4) we arrive at (2.5) and (2.6) after simple computation.

### 3. The Todorov - Charalambides identity

Here we shall discuss equation (1.4). Let $s(n,k)$ and $S(n,k)$ be the Stirling numbers of the first kind and the second kind correspondingly (see [9]). It is known that these numbers satisfy the orthogonality relation

$$\sum_{k=0}^n s(n,k)\, S(k,m) = \delta_{mn} = \begin{cases} 0 & m \neq n \\ 1 & m = n \end{cases}$$

while the alternating sums are related to the Lah numbers (see [9], p.156):

(3.1) $$L(n,m) = (-1)^n \sum_{k=0}^{n} s(n,k) S(k,m)(-1)^k.$$

The following identity extends this representation.

**Proposition 3**. For any two nonnegative integers $n, m$, and every complex number $z$ we have

(3.2) $$\frac{m!}{n!} \sum_{k=0}^{n} s(n,k) S(k,m) z^k = (-1)^m \sum_{j=0}^{m} \binom{m}{j} (-1)^j \binom{zj}{n}.$$

This identity was obtained by Todorov [19], who showed that both sides equal Taylor's coefficients of the function $f(t) = ((1+t)^z - 1)^m$. It was also found independently by Charalambides in his study of the generalized factorial coefficients (see [7] and [8]). A short proof of (3.2) is given in the recent paper [3].

Now we show that equation (3.1) follows from (3.2). Setting $z = -1$ in (3.2) we find

(3.3) $$\frac{m!}{n!} \sum_{k=0}^{n} s(n,k) S(k,m)(-1)^k = (-1)^m \sum_{j=0}^{m} \binom{m}{j} (-1)^j \binom{-j}{n}.$$

The right hand side can be written in the form

$$(-1)^{m+n} \sum_{j=0}^{m} \binom{m}{j} (-1)^j \binom{n+j-1}{n}$$

since

$$\binom{-j}{n} = (-1)^n \binom{n+j-1}{n}.$$

Next we use a well-known identity from [11]

(3.4) $$\sum_{j=0}^{m} \binom{m}{j} (-1)^j \binom{y+j}{n} = (-1)^m \binom{y}{n-m}$$

and choosing $y = n - 1$ we find

(3.5) $$(-1)^{m+n}\sum_{j=0}^{m}\binom{m}{j}(-1)^j\binom{n+j-1}{n}=(-1)^n\binom{n-1}{n-m}=(-1)^n\binom{n-1}{m-1}.$$

Now from (3.3)

$$\frac{m!}{n!}\sum_{k=0}^{n}s(n,k)\,S(k,m)(-1)^k=(-1)^n\binom{n-1}{m-1},$$

or

$$\sum_{k=0}^{n}s(n,k)\,S(k,m)(-1)^k=(-1)^n L(n,m)$$

which proves (3.1).

At the same time we can apply (3.5) to equation (1.5). This gives

$$D^n e^{1/x}=(-1)^n n!\,e^{1/x}x^{-n}\sum_{k=1}^{n}\frac{(-1)^k}{k!}x^{-k}\left\{\sum_{j=1}^{k}(-1)^j\binom{k}{j}\binom{n+j-1}{n}\right\}$$

$$=(-1)^n n!\,e^{1/x}x^{-n}\sum_{k=1}^{n}\frac{(-1)^k}{k!}x^{-k}\left\{(-1)^k\frac{k!}{n!}L(n,k)\right\}$$

$$=(-1)^n e^{1/x}x^{-n}\sum_{k=1}^{n}L(n,k)x^{-k}$$

which is exactly (1.2). We see that the formula for the derivatives $D^n e^{1/x}$ was practically found by Schwatt.

## 4. Schwatt's formula in terms of exponential polynomials

With the help of the Todorov - Charalambides identity, Schwatt's formula (1.4) can be written in terma of Stirling numbers and exponential polynomials.

The polynomials $\varphi_n(x), n=0,1,...,$ defined by

$$\varphi_m(x) = \sum_{k=0}^{m} S(m,k)\, x^k$$

are known as the exponential polynomials (or single-variable Bell polynomials) – see [4] and [5]. They have the generating function

$$e^{x(e^t-1)} = \sum_{n=0}^{\infty} \varphi_n(x) \frac{t^n}{n!} ,$$

and can be defined also by the important property $(xD)^n e^x = \varphi_n(x)\, e^x$, $n = 0,1,\ldots$.

**Proposition 4.** For any $n \geq 0$ and any two numbers $c, p$ we have

(4.1) $\quad D^n e^{c x^p} = e^{c x^p} x^{-n} \sum_{j=0}^{n} s(n,j)\, p^j \varphi_j(c\, x^p)$

and in particular, when $c = 1$ and $p = -1$,

(4.2) $\quad D^n e^{1/x} = e^{1/x} x^{-n} \sum_{j=0}^{n} s(n,j)(-1)^j \varphi_j(1/x)$ .

*Proof.* Substituting (3.2) in (1.4) we obtain

$$D^n e^{c x^p} = n!\, e^{c x^p} x^{-n} \sum_{k=1}^{n} \frac{(-1)^k}{k!} c^k x^{pk} \left\{ (-1)^k \frac{k!}{n!} \sum_{j=0}^{k} s(n,j)\, S(j,k)\, p^j \right\}$$

$$= e^{c x^p} x^{-n} \sum_{k=1}^{n} c^k x^{pk} \left\{ \sum_{j=0}^{k} s(n,j)\, S(j,k)\, p^j \right\}$$

$$= e^{c x^p} x^{-n} \sum_{j=0}^{n} s(n,j)\, p^j \left\{ \sum_{k=1}^{j} c^k x^{pk} S(j,k) \right\}$$

$$= e^{c x^p} x^{-n} \sum_{j=0}^{n} s(n,j)\, p^j \varphi_j(c\, x^p) .$$

Comparing this to (1.2) we arrive at the identity

$$(4.3) \quad \sum_{k=1}^{n} L(n,k) x^k = (-1)^n \sum_{j=0}^{n} s(n,j)(-1)^j \varphi_j(x).$$

Also, from (2.3),

$$(4.4) \quad L_n^{(-1)}(x) = \frac{(-1)^n}{n!} \sum_{j=0}^{n} s(n,j)(-1)^j \varphi_j(-x).$$

With $p = 1/2$ in (4.1) we have

$$(4.5) \quad D^n e^{c\sqrt{x}} = e^{c\sqrt{x}} x^{-n} \sum_{j=0}^{n} \frac{1}{2^j} s(n,j) \varphi_j(c\sqrt{x}).$$